\documentclass[12pt]{amsart}
\usepackage{amssymb,amsmath,latexsym}
\theoremstyle{plain}
\newtheorem{theorem}{Theorem}

\newtheorem{corollary}[theorem]{Corollary}

\newtheorem{lemma}[theorem]{Lemma}
\theoremstyle{definition}
\newtheorem{definition}[theorem]{Definition}
\newtheorem{Example}[theorem]{Example}

\newtheorem*{Notation}{Notation}

\begin{document}
\title[Kronecker Product]{The Kronecker Product of Schur 
Functions indexed by Two-Row Shapes or  Hook Shapes.}
\author{Mercedes H. Rosas}
\address{       Department of Mathematics,
        	 Brandeis University\\
               Waltham, MA 02254}
\email{rosas@math.brandeis.edu}
\begin{abstract}

The Kronecker product of two Schur functions $s_{\mu}$ and 
$s_{\nu}$, denoted by $s_{\mu}*s_{\nu}$, is the Frobenius 
characteristic of the tensor product of the irreducible 
representations of the symmetric group
corresponding to the partitions $\mu$ 
and $\nu$. The coefficient of $s_{\lambda}$ in this
product is denoted by $\gamma^{\lambda}_{{\mu}{\nu}}$,
and corresponds to the multiplicity of the irreducible 
character $\chi^{\lambda}$ in $\chi^{\mu}\chi^{\nu}.$

We use Sergeev's Formula for a Schur
function of a difference of two alphabets and the
comultiplication expansion for $s_{\lambda}[XY]$ to find
closed formulas for the Kronecker coefficients
$\gamma^{\lambda}_{{\mu}{\nu}}$ when
$\lambda$ is an arbitrary shape and $\mu$ and $\nu$ are 
 hook shapes  or  two-row shapes.

Remmel \cite{Re1, Re2} and Remmel and Whitehead \cite{Re-Wh} 
 derived some closed formulas for the Kronecker
product of Schur functions indexed by two-row shapes or hook shapes using
a different approach.  We believe that the approach of this paper is more 
natural. The formulas obtained are simpler and reflect the symmetry of 
the Kronecker product.
\end{abstract}
\maketitle


\catcode `!=11

\newdimen\squaresize 
\newdimen\thickness 
\newdimen\Thickness
\newdimen\ll! \newdimen \uu! \newdimen\dd! \newdimen \rr! \newdimen
\temp!

\def\sq!#1#2#3#4#5{%
\ll!=#1 \uu!=#2 \dd!=#3 \rr!=#4
\setbox0=\hbox{%
 \temp!=\squaresize\advance\temp! by .5\uu!
 \rlap{\kern -.5\ll! 
 \vbox{\hrule height \temp! width#1 depth .5\dd!}}%
%
 \temp!=\squaresize\advance\temp! by -.5\uu!  
 \rlap{\raise\temp! 
 \vbox{\hrule height #2 width \squaresize}}%
%
 \rlap{\raise -.5\dd!
 \vbox{\hrule height #3 width \squaresize}}%
%
 \temp!=\squaresize\advance\temp! by .5\uu!
 \rlap{\kern \squaresize \kern-.5\rr! 
 \vbox{\hrule height \temp! width#4 depth .5\dd!}}%
%
 \rlap{\kern .5\squaresize\raise .5\squaresize
 \vbox to 0pt{\vss\hbox to 0pt{\hss $#5$\hss}\vss}}%
}
 \ht0=0pt \dp0=0pt \box0
}

\def\vsq!#1#2#3#4#5\endvsq!{\vbox to \squaresize{\hrule
width\squaresize height 0pt%

\vss\sq!{#1}{#2}{#3}{#4}{#5}}}

\newdimen \LL! \newdimen \UU! \newdimen \DD! \newdimen \RR!

\def\vvsq!{\futurelet\next\vvvsq!}
\def\vvvsq!{\relax
  \ifx     \next l\LL!=\Thickness \let\continue!=\skipnexttoken!
  \else\ifx\next u\UU!=\Thickness \let\continue!=\skipnexttoken!
  \else\ifx\next d\DD!=\Thickness \let\continue!=\skipnexttoken!
  \else\ifx\next r\RR!=\Thickness \let\continue!=\skipnexttoken!
  \else\ifx\next P\let\continue!=\place!
  \else\def\continue!{\vsq!\LL!\UU!\DD!\RR!}%
  \fi\fi\fi\fi\fi 
  \continue!}

\def\skipnexttoken!#1{\vvsq!}

\def\place! P#1#2#3{%
\rlap{\kern.5\squaresize\temp!=.5\squaresize\kern#1\temp!
  \temp!=\squaresize \advance\temp! by #2\squaresize \temp!=.5\temp!
  \raise\temp!\vbox to 0pt{\vss\hbox to 0pt{\hss$#3$\hss}\vss}}\vvsq!}

\def\Young#1{\LL!=\thickness \UU!=\thickness \DD! = \thickness \RR! =
\thickness
\vbox{\smallskip\offinterlineskip
\halign{&\vvsq! ## \endvsq!\cr #1}}}

\def\blank{\omit\hskip\squaresize}
\catcode `!=12


\def\picture #1 by #2 (#3){
  \vbox to #2{
    \hrule width #1 height 0pt depth 0pt
    \vfill
    \special{picture #3} 
    }
  }

\def\scaledpicture #1 by #2 (#3 scaled #4){{
  \dimen0=#1 \dimen1=#2
  \divide\dimen0 by 1000 \multiply\dimen0 by #4
  \divide\dimen1 by 1000 \multiply\dimen1 by #4
  \picture \dimen0 by \dimen1 (#3 scaled #4)}
  }

					\section{Introduction}

The aim of this paper is to derive an explicit formula for the
Kronecker coefficients corresponding to partitions of certain shapes.
The Kronecker coefficients, $\gamma^{\lambda}_{{\mu}{\nu}}$, arise
when  expressing a Kronecker product (also called inner or internal
product), 
$s_{\mu}*s_{\nu}$, of Schur functions in the Schur basis,
	\begin{equation}
		s_{\mu}*s_{\nu} = 
		\sum_{\mu,\nu} \gamma^{\lambda}_{{\mu}{\nu}}s_{\lambda}.
		\label{eq:first}
	\end{equation}
These coefficients can also be defined as the multiplicities of the 
irreducible representations in the tensor product of two irreducible 
representations of the symmetric group. A  third way to define them 
is by the comultiplication expansion. Given two alphabets
 $X=\{x_1, x_2, \cdots\}$ and $Y=\{y_1, y_2, \cdots\}$
	\begin{equation}
		s_{\lambda}[XY] = 
		\sum_{\mu,\nu} \gamma^{\lambda}_{{\mu}{\nu}}
		s_{\mu}[X]s_{\nu}[Y],
		\label{eq:main}
	\end{equation}
where $s_{\lambda}[XY]$ means 
$s_{\lambda}(x_1y_1, x_1y_2,\cdots,x_iy_j, \cdots)$.
Remmel \cite{Re1, Re2} and Remmel and Whitehead \cite{Re-Wh}  
have studied the Kronecker product of Schur functions corresponding 
to two two-row shapes, two hook shapes, and a hook shape and a two-row 
shape. We will use the comultiplication expansion (\ref{eq:main})
for the Kronecker coefficients, and a formula for expanding a Schur
 function of a difference of two alphabets due to Sergeev \cite{Be-Ga}
 to obtain similar results in a simpler way. We believe that the 
formulas obtained using this approach are elegant and reflect the symmetry of
the Kronecker product.
In the three cases
we found a way to express the Kronecker coefficients in terms of regions
and paths in  ${\bf N}^2$.

						\section{Basic definitions} 

A partition $\lambda$ of a positive integer $n$, written 
as $\lambda \vdash n$, is an unordered sequence  of natural
 numbers adding to $n$.  We write $\lambda$ as
 $ \lambda = (\lambda_1,\lambda_2,\cdots,\lambda_n)$, 
where  $\lambda_1 \ge \lambda_2 \ge \cdots$, and consider
two  such strings equal if they differ by a string of
zeroes.   The nonzero numbers $\lambda_i$ are  called the
parts of 
$\lambda$, and the number of parts is called  the length 
of $\lambda$, denoted by   $l(\lambda)$. 
In some cases, it is  convenient to write 
$ \lambda =(1^{d_1} 2^{d_2} \cdots n^{d_n})$ 
for the partition  of $n$ that has $d_i$ copies of $i$. 
Using this notation, we  define the integer $z_{\lambda}$ to
 be $1^{d_1}d_1!\, 2^{d_2}d_2!\,  \cdots n^{d_n}d_n!$.

We identify  $\lambda$ with the set of
 points $(i,j)$ in ${\bf N}^2$ defined by
$1 \le j \le \lambda_i$, and refer to them as the Young 
diagram of $\lambda$. The Young diagram of a partition
$\lambda$ is thought of as a collection of boxes  
arranged using matrix coordinates.
 For instance, the Young diagram corrresponding to 
$\lambda=(4,3,1)$ is 
$$
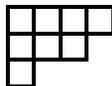
\squaresize = 10pt
\thickness = 1pt
\Thickness = .2pt
\Young{
& & &   \cr
& &   \cr
  \cr
}$$
To any partition $\lambda$ we associate the partition 
$\lambda'$, its conjugate
partition, defined by $\lambda_i'= \left| \{j : \lambda_j \ge i \}
\right|.$ Geometrically, $\lambda'$ can be obtained from
 $\lambda$  by flipping the Young diagram  of $\lambda$ around its
 main diagonal. For instance, the conjugate partition of $\lambda$ is 
$\lambda'=(3,2,2,1)$, and the corresponding Young diagram is
$$
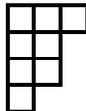
\squaresize = 10pt
\thickness = 1pt
\Thickness = .2pt
\Young{
& &  \cr
&    \cr
&  \cr
  \cr
}$$

We recall some facts about the theory of representations
 of the symmetric group, and about symmetric functions.
 See \cite{Mac} or \cite{Sag} for proofs and details.

Let $R(S_n)$ be the space of class function in $S_n$, 
the  symmetric group on $n$ letters, and let $\Lambda^n$
 be the space of homogeneous symmetric  functions of 
degree $n$. A  basis for $R(S_n)$ is given by the characters 
of the irreducible representations of
$S_n$. Let $\chi^{\mu}$ be the irreducible character 
of $S_n$ corresponding to the partition
$\mu$. There is a scalar product 
$\langle \, \, , \, \rangle_{S_n}$ on $R(S_n)$ defined 
by
\[
\langle \chi^{\mu},\chi^{\nu} \rangle_{S_n} = 
\frac{1}{n!}\sum_{\sigma \in S_n}
		 \chi^{\mu}(\sigma) \chi^{\nu}(\sigma),
\]
 and extended by
linearity.  

A basis for the space of symmetric
functions is given by the Schur functions.
There exists a scalar product
 $\langle \, \, , \, \rangle_{\Lambda^n}$ on $\Lambda^n$ 
defined by 
\[
\langle s_{\lambda},s_{\mu} 
\rangle_{\Lambda^n}=\delta_{\lambda \mu},
\]
where $\delta_{\lambda \mu}$ is the Kronecker delta, and  extended by
linearity.  

Let $p_{\mu}$ be the power sum symmetric function 
corresponding to $\mu$, where $\mu$ is a partition of $n$. 
There is an isometry
 $\mathop{\rm ch}^n:R(S_n)
\mapsto \Lambda^n$,
 given by the characteristic map,
\[
			{\mathop{\rm ch}}^n (\chi)=
			\sum_{\mu \vdash n} z_{\mu}^{-1} 
			\chi(\mu)p_{\mu}.
 \]
This map has the remarkable property that if $\chi^{\lambda}$ is the
 irreducible character of $S_n$ indexed by $\lambda$, then
 $\mathop{\rm ch}^n(\chi^{\lambda})=s_{\lambda}$, the Schur function 
corresponding to $\lambda$. In particular, we obtain that 
$s_{\lambda}=
			\sum_{\mu \vdash n} z_{\mu}^{-1} 
			\chi^{\lambda}(\mu)p_{\mu}.$ Hence, 
\begin{equation}
\chi^{\lambda}(\mu)=\langle  s_{\lambda}, p_{\mu} 
\rangle.
\label{eq:character}
\end{equation}
Let $\lambda$, $\mu$, and $\nu$ be partitions of $n$. The Kronecker 
 coefficients $\gamma^{\lambda}_{{\mu}{\nu}}$ are defined by
\begin{align}
 	\gamma^{\lambda}_{{\mu}{\nu}} =
	\langle \chi^{\lambda}, \chi^{\mu} \chi^{\nu} {\rangle}_{S_n}
	= \frac{1}{n!} 
		\sum_{\sigma \in S_n}
		\chi^{\lambda}(\sigma) \chi^{\mu}(\sigma)
		\chi^{\nu}(\sigma).\label{eq:defikron}
\end{align}
Equation (\ref{eq:defikron}) shows that the Kronecker coefficients
$
\gamma^{\lambda}_{{\mu}{\nu}}$ are symmetric in
 $	\lambda$, $\mu$, and $\nu$.
The relevance of the Kronecker coefficients comes from the
following fact: Let $X^{\mu}$ be the representation of the 
symmetric group corresponding
to the character
 $\chi^{\mu}$. Then $\chi^{\mu}\chi^{\nu}$ is  
 the character of $X^{\mu} \otimes X^{\nu}$, the representation
 obtained by taking the tensor product of $X^{\mu}$ and $X^{\nu}$. 
Moreover, $\gamma^{\lambda}_{{\mu}{\nu}}$ is the multiplicity 
of $X^{\lambda}$ in $X^{\mu} \otimes X^{\nu}$.

Let $f$ and $g$ be homogeneous symmetric functions of degree $n$. The
Kronecker product,
 $f*g$, is defined by
		\begin{equation}
			f*g={\mathop{\rm ch}}^n(uv),
		\label{eq:star}
		\end{equation}
where $\mathop{\rm ch}^nu=f$, $\mathop{\rm ch}^nv=g$, and 
$uv(\sigma)=u(\sigma)v(\sigma)$. To obtain (\ref{eq:first}) 
from this definition, we set $f=s_{\mu}$,
 $g=s_{\nu}$, $u=\chi^{\mu}$, and $v=\chi^{\nu}$ in (\ref{eq:star}).
 
The Kronecker product has the following symmetries:
\begin{align*}
s_{\mu}*s_{\nu}&=s_{\nu}*s_{\mu}.\\
s_{\mu}*s_{\nu}&=s_{\mu'}*s_{\nu'}.
\end{align*}
Moreover,  if $\lambda$ is a one-row shape
\[
\gamma^{\lambda}_{{\mu}{\nu}}=\delta_{\mu,\nu}.
\]    
We introduce the operation of substitution or 
plethysm into a symmetric function. 
Let $f$ be a symmetric function, and let 
$X=\{x_1,x_2,\cdots\}$ be an alphabet.
 We write $X=x_1+x_2+\cdots$, and define $f[X]$ by,
\[
f[X]=f( x_1,x_2,\cdots).
\nonumber
\]
In general, if $u$ is any element of
${\bf Q}[[x_1, x_2,
\cdots]]$, we write $u$ as 
$\sum_{\alpha}{c_{\alpha}u_{\alpha}}$ where $u_{\alpha}$ is a
monomial with coefficient $1$. Then $p_{\lambda}[u]$ is defined 
by setting
\begin{align}
p_n[u]&=\sum_{\alpha}{c_{\alpha}u_{\alpha}^n}
\nonumber\\
p_{\lambda}[u]&=p_{\lambda_1}[u]
\cdots p_{\lambda_n}[u]
\nonumber
\end{align}
 for $\lambda=(\lambda_1, \cdots, \lambda_n)$. We define $f[u]$ for
all symmetric functions $f$ by saying that
$f[u]$ is linear in $f$. 

The operation of substitution into a symmetric
function has the following properties.
For $\alpha$ and $\beta$ rational numbers, $ 
(\alpha f+\beta g)[u]=\alpha f[u]+\beta
g[u].$ Moreover, if $c_{\alpha}=1$ for all $\alpha$, then $
f[u]=f(\cdots,u_{\alpha},\cdots).$

Let $X=x_1+x_2+\cdots$ and $Y=y_1+y_2+\cdots$ be two
alphabets.   Define the sum of two alphabets by
$X+Y= x_1+ x_2+ \cdots+ y_1+ y_2+ \cdots,$
 and the product of two alphabets by  
$XY=x_1y_1+\cdots+x_iy_j+\cdots .$ 
Then
		\begin{align}
			p_n[X+Y] &= p_n[X]+p_n[Y],\nonumber\\
			p_n[XY] &=p_n[X]p_n[Y].
		\label{eq:multiplicative}
		\end{align} 
The inner product of function in the space of
symmetric functions in two infinite alphabets is defined by
\[
\langle \;\; , \; {\rangle}_{_{XY}}=\langle \;\; , \; {\rangle}_{_{X}}
\langle \;\; , \; {\rangle}_{_{Y}},
\]
where for any given alphabet $Z$,
$\langle \;\; , \; {\rangle}_{_{Z}}$ denotes the inner product of the 
space of symmetric functions in $Z$.

For all partitions $\rho$, we have that $p_{\rho}[XY]=p_{\rho}[X]p_{\rho}[Y].$
If we rewrite (\ref{eq:character})  as
$p_{\rho}=\sum_{\lambda}{\chi^{\lambda}(\rho)}s_{\lambda}$, then
\begin{equation}
\sum_{\lambda}{\chi^{\lambda}s_{\lambda}[XY]}=
\sum_{\mu,\nu} \chi^{\mu}\chi^{\nu}s_{\mu}[X]s_{\nu}[Y].
\label{eq:equiv1}
\end{equation}
Taking the coefficient of  ${\chi}^{\lambda}$ on 
both sides of the previous  equation
 we obtain
$$s_{\lambda}[XY]=\sum{\langle{\chi}^{\lambda},
{\chi}^{\mu}{\chi}^{\nu}\rangle s_{\mu}[X]s_{\nu}[Y]}.$$
Finally, using the definition of Kronecker coefficients (\ref{eq:defikron})
we obtain the comultiplication expansion (\ref{eq:main}).

\begin{Notation} 
Let $p$ be a point in ${\bf N}^2.$ We say  that $(i,j)$ 
can be reached from $p$, written
$p \leadsto (i,j)$, if $(i,j)$ can be reached from $p$ by moving any
number of steps south-west or north-west.
We define the weight function $\omega$ by
\[
\omega_p(i,j)=
\begin{cases}
x^iy^j, &\text{ if } p \leadsto (i,j),\\
0, &\text{otherwise.}
\end{cases}
\]
In particular, $\sigma_{k,l}(h)=0$ if $h<0$.
\end{Notation}

\begin{Notation}  

We denote by $\lfloor x \rfloor$ the largest integer less than
or equal to $x$ and by $\lceil x \rceil$ the smallest integer
greater than or equal to $x$.

If $f$ is a formal power series, then $[x^{\alpha}]\,f$ denotes
the coefficient of $x^{\alpha}$ in $f.$

Following Donald Knuth we denote the characteristic function applied to
a  proposition $P$ by enclosing $P$ with brackets,
\begin{equation}
(P)=
\begin{cases}
1, &\text{if proposition $P$ is true,}\\
0, &\text{otherwise.}
\end{cases}
\nonumber
\end{equation}
\end{Notation}

\section{The case of two two-row shapes}

The object of this section is to find a closed formula for
 the Kronecker coefficients  when
$\mu=(\mu_1, \mu_2)$ and  $\nu=(\nu_1, \nu_2)$ are 
two-row shapes, and when we do not have any restriction on the
partition $\lambda$. We describe the Kronecker coefficients 
$\gamma^{\lambda}_{{\mu}{\nu}}$ in terms of paths in ${\bf N}^2$.
More precisely, we define two rectangular
regions in ${\bf N}^2$ using the parts of $\lambda$.
Then we count the number of points in ${\bf N}^2$ inside each of these 
rectangles that can be reached from $(\nu_2, \mu_2+1)$, 
if we are allowed to move any number of steps south-west
or north-west.
Finally, we subtract these two numbers.

We begin by introducing two lemmas that 
allow us to state Theorem \ref{le:2rowst} in a concise 
form.

\begin{Notation} 
We use the coordinate axes as if we were working with matrices with first
entry $(0,0)$. That is, the point $(i,j)$ belongs to the $i$th row and 
the $j$th column.
\end{Notation}

\begin{lemma}\label{sigma}
Let $k$ and $l$ be positive numbers. Let $R$ be the rectangle with
width $k$, height $l$, and upper--left square $(0,0)$. Define
\[
\sigma_{k,l}(h)=\big | \{(u,v) \in R \cap {\bf N}^2 : (0,h)\leadsto (u,v) \} 
\big |
\]
Then
\begin{equation}
\sigma_{k,l}(h)=
\begin{cases}
0, &\text{if } h <0\\
\lfloor (\frac{h}{2}+1)^2 \rfloor, 
	&\text{if } 0\le h<\min (k,l)\\
\sigma_{k,l}(s)+(\frac{h-s}{2})\min(k,l), 
	&\text{if  } \min (k,l)\le h<\max (k,l)\\
\lceil \frac{kl}{2}\rceil - \sigma_{k,l}(k+l-h-4),
	&\text{if $h$ is even and } \max (k,l)\le h\\
\lfloor \frac{kl}{2} \rfloor - \sigma_{k,l}(k+l-h-4),
	&\text{if $h$ is odd and } \max (k,l)\le h\\
\end{cases}
\nonumber
\end{equation}
where $s$ is defined as follows: If $h - \min(k,l)$ is even, then $s=\min(k,l)-2$;
otherwise $s=\min(k,l)-1$.
\end{lemma}

\begin{proof}
If $h$ is to the left of the $0$th column, then we cannot reach any of the
points in ${\bf N}^2$ inside $R$. Hence, $\sigma_{k,l}(h)$ should be equal to zero.

If $0 \le h \le \min(k,l)$, then we are counting the number of points in ${\bf N}^2$
that can be reached from $(0,h)$ inside the square $S$ of side $\min(k,l)$.
We have to consider two cases. 
If $h$ is odd, then we are summing
$2+4+\cdots+(h+1)=\lfloor (\frac{h}{2}+1)^2\rfloor$.
On the other hand, if $h$ is even, then we are summing 
$1+3+\cdots+(h+1)=(\frac{h}{2}+1)^2$.

If $\min(k,l)\le h <\max(k,l)$, then we subdivide our problem into two parts.
First, we count the number of points in ${\bf N}^2$ that can be reached from $(0,h)$ 
inside the square $S$ by $\sigma_{k,l}(s)$.
Then we count those points in ${\bf N}^2$ that are in $R$ but not in $S$. Since 
$h<\max(k,l)$ all diagonals have length $\min(k,l)$ and there
are $\frac{h-s}{2}$ of them. See Table 1.

If $\max(k,l)\le h$, then it is easier to count the total number of points in ${\bf N}^2$
that can be reached from $(0,h)$ inside $R$  by choosing another parameter 
$\hat h $  big enough and with the  same parity as $h$. 
Then we subtract those points in ${\bf N}^2$ in $R$ that are not reachable from
$(0,h)$ because $h$ is too close.
If $h$ is even this number is $\lceil kl/2 \rceil$. If $h$ is 
odd this number is  $\lfloor kl/2 \rfloor$.

Then we subtract those points that we should not have counted.
we express this number in terms of the function $\sigma$.
The line $y=-x+h+2$ intersects the line $y=l-1$ 
at $x=h-l+3$. This is the $x$ coordinate of the first point on the 
last row that is not reachable from $(0,h)$. Then to obtain the number 
of points that can be reached from this point by moving south-west or
north-west, we subtract $h-l+3$ to $k-1$.  We have obtained that are
$\sigma_{k,l}(k+l-h-4)$ points that we should not have counted.
\end{proof}

\begin{Example} 

By definition $\sigma_{9,5}(4)$ counts the points in ${\bf N}^2$ in Table 1 marked with 
 $\circ$. Then $\sigma_{9,5}(4)=9$. Similarly,  $\sigma_{9,5}(8)$
counts the points in ${\bf N}^2$ in Table 1 marked either with the symbol $\circ$ 
or with the symbol $\bullet$.
Then  $\sigma_{9,5}(8)=19.$ 

\medskip

\centerline{
\begin{tabular}{|c|c|c|c|c|c|c|c|c|}
\hline
$\circ$ & & $\circ$ &  &  $\circ$
& \phantom{$\circ$} & $\bullet$ &  & $\bullet$\\
\hline
& $\circ$ & & $\circ$ & \phantom{$\circ$}
 & $\bullet$ & &$\bullet$ & \\ 
\hline 
$\circ$ &  & $\circ$ &  & $\bullet$
 & & $\bullet$ & & \\ 
 \hline
 & $\circ$ &  &$\bullet$ & 
 & $\bullet$ &  & & \\ 
 \hline
$\circ$ &  & $\bullet$ &  & $\bullet$
 & \phantom{$\circ$} & & & \\ \hline
\end{tabular}
}
\smallskip
{\centerline{Table 1. }} 

\end{Example}

\begin{lemma}\label{le:2rows2}
Let $a, b, c,$ and $d$ be in $\bf N$. Let $R$ be the 
rectangle with vertices $(a,c)$,
$(a+b,c)$, $(a,c+d)$, and $(a+b,c+d)$. We define
\[
\Gamma(a,b,c,d)(x,y)=\big | \{(u,v) \in R : (x,y)\leadsto (u,v) \} 
\big |.
\]
Suppose that $(x,y)$ is such that $x \ge y$.
Then 
\[
\Gamma(a,b,c,d)(x,y)=
\begin{cases}
\sigma_{b+1,d+1}(x+y-a-c), & 0 \le y \le c\\
\sigma_{b+1,y-c+1}(x-a)+ \sigma_{b+1,c+d-y+1}(x-a)
-\delta, & c< y < c+d\\
\sigma_{b+1,d+1}(x-y+c+d-a), & c+d \le y
 \end{cases}
\]
where $\delta$ is defined as follows
If $x<a$, then $\delta=0$.
If $a \le x \le a+b$, then $\delta=\big\lceil \frac{x-a+1}{2} \big\rceil$.
Finally, if $x>a+b$ then we consider two cases: 
If $x-a-b$ is even then $\delta=\big\lceil \frac{b+1}{2} \big\rceil$; otherwise,
 $\delta=\big\lfloor \frac{b+1}{2} \big\rfloor$.
\end{lemma}
\begin{proof}
We consider three cases. If $0\le y\le c$ then the first position inside
$R$ that we reach is $(x+y-a-c,c)$. Therefore, we assume that we are starting 
at this point. Similarly, if $y\ge c+d$, then the first position  inside
$R$ that we reach is $(x-y+c+d-a,c)$. Again, we can assume that we are starting 
at this point.

On the other hand, if $c<y<c+d$, then we subdivide the problem in two parts.
The number of position to the north of us is counted by $\sigma_{b+1,y-c+1}(x-a)$.
The number of position to the south of us is counted by $\sigma_{b+1,c+d-y+1}(x-a)$.
We define $\delta$ to be the number of points in ${\bf N}^2$ that we counted twice during
this process. Then it is easy to see that $\delta$ is given by the 
previous definition.
\end{proof}

To compute the coefficient
 $u_{\nu}$ in the expansion
$f[X]=\sum_{\eta}{u_{\eta}s_{\eta}[X]}$ for $f \in
{\Lambda}$, it is enough to expand $f[x_1+\cdots+x_n]=
	\sum_{\eta} {u_{\eta}s_{\eta}[x_1+\cdots+x_n]}$
for any $n \geq l({\nu})$. (See [7, section I.3],
 for proofs and details.)  
Therefore, in this section we work with symmetric functions
 in a finite number of variables.

Let $\mu$ and $\nu$ be two-row partitions. Set   $X=1+x$ and
$Y=1+y$ in the comultiplication expansion
(\ref{eq:main}) to obtain
	\begin{equation}
		s_{\lambda}[(1+y)(1+x)] = 
		\sum_{\mu,\nu} \gamma^{\lambda}_{{\mu}{\nu}}
		s_{\mu}[1+y]s_{\nu}[1+x].
		\label{eq:main2row}
	\end{equation}
Note that the Kronecker coefficients are zero when
$l(\lambda)>4.$ 

Jacobi's  definition of a Schur  function on a finite alphabet
$s_\lambda[X]$ as a quotient  of alternants says that   
 	\begin{align}
		s_{\lambda}[X]=
  		s_{\lambda}(x_1, \cdots, x_n)
	=\frac{\mathop{\rm det}
  		(x_i^{\lambda_j+n-j})_{1 \leq i,j\leq n}}
		{\prod_{i<j}(x_i-x_j)}.
	\label{eq:alternant}
	\end{align}

By the symmetry properties of the Kronecker product  it  is
 enough to compute the Kronecker coefficients
$\gamma_{\mu \nu}^{\lambda}$ when $ \nu_2 \le \mu_2$.

\begin{theorem}\label{le:2rowst}
Let $\mu$, $\nu$, and $\lambda$ be partitions of $n$, 
where $\mu=(\mu_1,\mu_2)$ and  
$\nu=(\nu_1, \nu_2)$ are two two-row partitions and let 
$\lambda=(\lambda_1,\lambda_2,\lambda_4,\lambda_4)$
be a partition of length less than or equal to $4$.
 Assume that $\nu_2 \le \mu_2 $. Then
\[
\gamma_{\mu\nu}^{\lambda}=\big(\Gamma(a,b,a+b+1,c)
                  -\Gamma(a,b,a+b+c+d+2,c)\big)(\nu_2,\mu_2+1).
\]
where $a=\lambda_3+\lambda_4$, 
$b=\lambda_2-\lambda_3$,
$c=\min(\lambda_1-\lambda_2, \lambda_3-\lambda_4)$ and 
$d=\big|\lambda_1+\lambda_4-\lambda_2-\lambda_3 \big |$.
		
\end{theorem}

\begin{proof} We expand the polynomial 
$
s_{\lambda}[(1+y)(1+x)]=s_{\lambda}(1,y,x,xy)
$
 in two different ways and obtain the Kronecker coefficients
by equating both results.
Let $\varphi$ be the polynomial defined by
 $\varphi=(1-x)(1-y)s_{\lambda}(1,y,x,xy).$
Using Jacobi's definition of a Schur 
function we obtain
\begin{equation}
  \varphi=
\frac{
		\begin{vmatrix}
   1&1&1&1\\
	y^{\lambda_1+3}&y^{\lambda_2+2}&y^{\lambda_3+1}
&y^{\lambda_4}\\
   x^{\lambda_1+3}&x^{\lambda_2+2}&x^{\lambda_3+1}
&x^{\lambda_4}\\
	(xy)^{\lambda_1+3}&(xy)^{\lambda_2+2}&(xy)^{\lambda_3+1}
&(xy)^{\lambda_4}
\end{vmatrix}  }
			{xy(1-xy)(y-x)(1-x)(1-y)}.
		\label{eq:2row1}
		\end{equation}
On the other hand, we may use Jacobi's definition to expand 
  $s_{\mu}[1+y]$ and $s_{\nu}[1+x]$ as  quotients of
 alternants. Substitute this into (\ref{eq:main2row}):
	\begin{align}
		s_{\lambda}[(1+y)(1+x)]&=
		\sum_{\substack{\mu=(\mu_1,\mu_2)\\ \nu=(\nu_1,\nu_2)}}
 {\gamma_{\mu \nu}^{\lambda}
  {\left( \frac{y^{\mu_2}-y^{\mu_1+1}} {1-y} \right) }
  {\left( \frac{x^{\nu_2}-x^{\nu_1+1}} {1-x} \right) }}
  \nonumber\\
 &=\sum_{\substack{\mu=(\mu_1,\mu_2)\\ \nu=(\nu_1,\nu_2)}}
		 \gamma_{\mu \nu}^{\lambda}
		\frac{
		x^{\nu_2}y^{\mu_2}-x^{\nu_2}y^{\mu_1+1}
		-x^{\nu_1+1}y^{\mu_2}+x^{\nu_1+1}y^{\mu_1+1}}
		{(1-x)(1-y)}.
		\label{eq:2row2}
	\end{align} 
Since $\nu_1+1$ and
$\mu_1+1$ are both greater than $ \lfloor \frac{n}{2} \rfloor,$ equation
 (\ref{eq:2row2}) implies that
{\em the coefficient of $x^{\nu_2}y^{\mu_2}$ in $\varphi$
 is $ \gamma_{\mu
\nu}^{\lambda}$. }
It is convenient to define an auxiliary polynomial by
\begin{equation}
	\zeta=(1-xy)(y-x) \varphi .
	 \label{eq:auxiliar}
\end{equation}
Let $\xi$ be the polynomial obtained by expanding the 
determinant appearing in  (\ref{eq:2row1}). Equations 
(\ref{eq:2row1}) and (\ref{eq:auxiliar}) imply
\[
 \zeta=\frac{\xi}{xy(1-x)(1-y)}.
\]

Let $\xi_{i,j}$ be the coefficient of $x^iy^j$ in
$\xi$. ( Then $\xi_{i,j}$ is zero if $i \le 0$
or $j \le 0$, because $\xi$ is a polynomial divisible by $xy$.)
Let $\zeta_{i,j}$ be the coefficient of
 $x^iy^j$ in $\zeta$. Then
	\begin{align}
		\sum_{i,j\ge0}{\zeta_{i,j}x^iy^j}=
		\frac{1}{xy(1-x)(1-y)}
		\sum_{i,j\ge0}{\xi_{i,j}x^i y^j}
		=\sum_{i,j,k,l\ge0}{\xi_{i-k,j-l}x^{i-1}y^{j-1}}.
		\label{eq:phi}
	\end{align}
Comparing the coefficient of $x^iy^j$ on both sides 
of equation (\ref{eq:phi}) we obtain that
\begin{align}
  \zeta_{i,j}=
  \sum_{k,l \geq 0}\xi_{i+1-k,j+1-l}
= \sum_{k=0}^{i}\sum_{l=0}^{j}{\xi_{k+1,l+1}}
  \label{eq:phi2}
\end{align}
We compute $\zeta_{i,j}$ from (\ref{eq:phi2}) by expanding the 
determinant appearing on (\ref{eq:2row1}).
We consider two cases.

{\bf Case 1.} Suppose that ${\lambda}_1+{\lambda}_4 
> {\lambda}_2 + {\lambda}_3$. 
Then
\[
	 {\lambda}_1+{\lambda}_2+4  >{\lambda}_1+
{\lambda}_3+3  > 
	{\lambda}_1+{\lambda}_4+2  \geq
	 {\lambda}_2+{\lambda}_3+2  > 
	{\lambda}_2+{\lambda}_4+1  >{\lambda}_3+{\lambda}_4.
\]
We record the values of  $\xi_{j+1, i+1}$ in 
 Table 2. We use the convention that $\xi_{i+1,j+1}$ 
is zero whenever the
$(i,j)$ entry is not in Table 2.

\bigskip

\centerline{
\begin{tabular}{ c c c c c c c }
\hline
$i \backslash j$ & \scriptsize${\lambda_3+\lambda_4}$ 
&\scriptsize${\lambda_2+\lambda_4+1}$ &
\scriptsize${\lambda_2+\lambda_3+2}$ &
\scriptsize${\lambda_1+\lambda_4+2}$ & 
\scriptsize${\lambda_1+\lambda_3+3}$ &
\scriptsize${\lambda_1+\lambda_2+4}$  \\ \hline
\scriptsize$\lambda_3+\lambda_4$ 
& \footnotesize{$\phantom{+}0$} & \footnotesize{$-1$} &
\scriptsize{$+1$} & \scriptsize{$+1$} 
& \scriptsize{$-1$} &
\scriptsize{$\phantom{+}0$} \\  
\scriptsize$\lambda_2+\lambda_4+1$ & 
\scriptsize{$+1$} &
\scriptsize{$\phantom{+}0$} &
\scriptsize{$-1$} & \scriptsize{$-1$} & 
\scriptsize{$\phantom{+}0$} &
\scriptsize{$+1$}  \\  
\scriptsize$\lambda_2+\lambda_3+2$ & 
\scriptsize{$-1$} &
\scriptsize{$+1$} & \scriptsize{$\phantom{+}0$} & 
\scriptsize{$\phantom{+}0$} &
\scriptsize{$+1$} &
\scriptsize{$-1$}  \\  
\scriptsize$\lambda_1+\lambda_4+2$ & 
\scriptsize{$-1$} &
\scriptsize{$+1$} & \scriptsize{$\phantom{+}0$} & 
\scriptsize{$\phantom{+}0$} &
\scriptsize{$+1$} &
\scriptsize{$-1$} \\  
\scriptsize$\lambda_1+\lambda_3+3$  & 
\scriptsize{$+1$} &
\scriptsize{$\phantom{+}0$} &
\scriptsize{$-1$} & \scriptsize{$-1$} & 
\scriptsize{$\phantom{+}0$} &
\scriptsize{$+1$}  \\  
\scriptsize$\lambda_1+\lambda_2+4$ & 
\scriptsize{$\phantom{+}0$} &
\scriptsize{$-1$} &
\scriptsize{$+1$} & \scriptsize{$+1$} & 
\scriptsize{$-1$} &
\scriptsize{$\phantom{+}0$} \\ \hline
\end{tabular}
}

\smallskip
{\centerline{Table 2}

\centerline{The values of $\xi_{j+1, i+1}$ when 
${\lambda}_1+{\lambda}_4 
\geq {\lambda}_2 + {\lambda}_3$}}

\medskip

Equation (\ref{eq:phi2}) shows that 
the value of $\zeta_{i,j}$ can be  obtained by adding 
the entries northwest of the point
$(i,j)$. In Table 3 we record the values of 
$\zeta_{i,j}$.

\medskip

\centerline{
\begin{tabular}{ c c c c c c c c }
\hline
$i \backslash j$& $I_1$ & $I_2$  & $I_3$  & $I_4$  & $I_5$ & $I_6$ &
$I_7$  \\ \hline
$I_1$ &$\phantom{+}0$ &$\phantom{+}0$ &$\phantom{+}0$
&$\phantom{+}0$ &$\phantom{+}0$ &$\phantom{+}0$ &
$\phantom{+}0$ \\
$I_2$ &$\phantom{+}0$ &$\phantom{+}0$ &$-1$ &$\phantom{+}0$
&$+1$ &$\phantom{+}0$ & $\phantom{+}0$ \\
$I_3$ &$\phantom{+}0$ &$+1$ &$\phantom{+}0$ &$\phantom{+}0$
&$\phantom{+}0$ &$-1$ & $\phantom{+}0$ \\
$I_4$ &$\phantom{+}0$ &$\phantom{+}0$ &$\phantom{+}0$
&$\phantom{+}0$ &$\phantom{+}0$ &$\phantom{+}0$ &
$\phantom{+}0$ \\
$I_5$ &$\phantom{+}0$ &$-1$ &$\phantom{+}0$ &$\phantom{+}0$
&$\phantom{+}0$ &$+1$ & $\phantom{+}0$ \\
$I_6$ &$\phantom{+}0$ &$\phantom{+}0$ &$+1$ &$\phantom{+}0$
&$-1$ &$\phantom{+}0$ & $\phantom{+}0$ \\
$I_7$ &$\phantom{+}0$ &$\phantom{+}0$ &$\phantom{+}0$
&$\phantom{+}0$ &$\phantom{+}0$ &$\phantom{+}0$ &
$\phantom{+}0$ \\
\hline
\end{tabular} }
\smallskip
{\centerline{Table 3}}
{\centerline{ The values of $\zeta_{i,j}$
when ${\lambda}_1+{\lambda}_4 
\geq {\lambda}_2 + {\lambda}_3$}
}

where
{\allowdisplaybreaks
\begin{align*}
I_1&=[0,\lambda_3+\lambda_4),
&I_2&=[\lambda_3+\lambda_4,\lambda_2+\lambda_4],\\
I_3&=[\lambda_2+\lambda_4+1,\lambda_2+\lambda_3+1],
&I_4&=[\lambda_2+\lambda_3+2,\lambda_1+\lambda_4+1],\\
I_5&=[\lambda_1+\lambda_4+2,\lambda_1+\lambda_3+2],
&I_6&=[\lambda_1+\lambda_3+3,\lambda_1+\lambda_2+3],\\ 
I_7&=[\lambda_1+\lambda_2+4, \infty). &\, &\,
\end{align*}}

{\bf Case 2.} Suppose that ${\lambda}_1+{\lambda}_4 \le
 {\lambda}_2 + {\lambda}_3$. 
Then
\[
	 {\lambda}_1+{\lambda}_2+4  >
{\lambda}_1+{\lambda}_3+3  > 
	{\lambda}_2+{\lambda}_3+2  >
 {\lambda}_1+{\lambda}_4+2  > 
	{\lambda}_2+{\lambda}_4+1  >{\lambda}_3+{\lambda}_4.
\]

Note that in Table 3, the rows and columns 
corresponding to ${\lambda}_1+{\lambda}_4 +2$ and
${\lambda}_2+{\lambda}_3+2$ are the same. Therefore, the
values of $\xi_{i,j}$ for ${\lambda}_1+{\lambda}_4 \le
 {\lambda}_2 + {\lambda}_3$
are recorded in Table 3, if we set 
\begin{align*}
&I_3=[{\lambda}_2+{\lambda}_4+1,{\lambda}_1+
{\lambda}_4+1]\\
&I_4=[{\lambda}_2+{\lambda}_4+2,{\lambda}_2+
{\lambda}_3+1]\\
&I_5=[{\lambda}_2+{\lambda}_3+2,{\lambda}_1+
{\lambda}_3+2],
\end{align*}
and define the other intervals as before. 

\medskip

In both cases, let $\varphi_{i,j}$ be the coefficient 
of $x^iy^j$ in $\varphi$. 
Using  (\ref{eq:auxiliar}) we obtain that
	\begin{align}
  	\varphi &=
		\frac{1}{(1-xy)(y-x)} 
			\sum_{i,j\ge0}{\zeta_{i,j}x^iy^j}
\nonumber\\
  		&=
  		\frac{1}{y-x}\sum_{i,j,l\ge0}
		{\zeta_{i-l,j-l}x^iy^j}
\label{eq:med}\\
		&=\sum_{i,j,k,l\ge0}{\zeta_{i-k-l,j+k-l+1}x^iy^j.}
\nonumber
	\end{align}
(Note: We can divide by $y-x$ because $\varphi=0$ when\, $x=y$.)
Comparing the coefficients of  $x^iy^j$ on both sides  
of equation (\ref{eq:med}), we  obtain 
$\varphi_{i,j}=\sum_{k,l\ge0}{\zeta_{i-k-l,j+k-l+1}}.$
Therefore,
\begin{equation}
\varphi_{\nu_2,\mu_2}=\sum_{i,j=0}^{\nu_2}
{\zeta_{\nu_2-i-j,\mu_2+i-j+1}}.
\label{eq:coef1}
\end{equation}
We have shown that $\gamma^{\lambda}_{\mu \nu}=\varphi_{\nu_2,\mu_2}$ 
can be obtained by adding the
entries in Table 3 in all points in ${\bf N}^2$ that can be reached from
$(\nu_2,\mu_2+1)$.
See Table 4.

\medskip

\centerline{
\begin{tabular}{|c|c|c|c|c|c|c|c|c|}
\hline
  &  &  & \phantom{$\circ$} &  \phantom{$\circ$} \\ 
  \hline
 $\circ$ &  &  & & \\
 \hline
   &$\circ$ &   & & \\
   \hline
 $\circ$ &  &$\circ$  & & \\ 
 \hline
  &$\circ$ &   & & \\
  \hline
 $\circ$ &    &  & & \\
 \hline
\end{tabular} }

\smallskip
{\centerline {Table 4}}
\medskip
The right-most point in Table 4 has coordinates $(3,2)$.

\medskip

By hypothesis $\nu_2 \le \mu_2 \le \lfloor n/2 \rfloor$. 
Then, if we start at $(\nu_2, \mu_2+1)$ and move as 
previously described, the only points in ${\bf N}^2$ that
we can possibly reach and that are
nonzero in Table 3 are those in $I_2\times I_3$ or $I_2\times I_5.$
Hence, we have  that  $\varphi_{\nu_2,\mu_2}$ 
is the number of points in ${\bf N}^2$
inside $I_2\times I_3$ that can be reached from $(\nu_2,\mu_2+1)$ minus
the ones that can be  reached  in  $I_2\times I_5.$ 

\medskip

{\bf Case 1} The inequality 
  ${\lambda}_1+{\lambda}_4 > {\lambda}_2 + {\lambda}_3$
implies that ${\lambda}_1+{\lambda}_4+1 > \lfloor \frac{n}{2}
\rfloor$. Moreover, 
$\mu_2 \ge \nu_2$ implies that we are only considering the region
of ${\bf N}^2$ given by $0 \leq i \leq j \leq  \lfloor \frac{n}{2}
\rfloor$. The number of points in ${\bf N}^2$ that can be reached 
from $(\nu_2,\mu_2+1)$ inside $I_2\times I_3$
is given by
$\Gamma(\lambda_3+\lambda_4,\lambda_2-\lambda_3,
\lambda_2+\lambda_4+1,\lambda_3-\lambda_4).$ Similarly, 
the number of points in ${\bf N}^2$ that can be reached 
from $(\nu_2,\mu_2+1)$ inside $I_2\times I_5$
is given by
$\Gamma(\lambda_3+\lambda_4,\lambda_2-\lambda_3,
\lambda_1+\lambda_4+2,\lambda_3-\lambda_4).$

\medskip

{\bf Case 2} The inequality
 $\lambda_2+\lambda_3 \ge \lambda_1+\lambda_4$,
 implies that ${\lambda}_1+{\lambda}_4+1 > \lfloor \frac{n}{2}
\rfloor$. Moreover, 
$\mu_2 \ge \nu_2$ implies that we are only considering the region
of ${\bf N}^2$ given by
$0 \leq i \leq j \leq  \lfloor \frac{n}{2}\rfloor$. 
The number of points in ${\bf N}^2$ that can be reached 
from $(\nu_2,\mu_2+1)$ inside $I_2\times I_3$
is given by
$\Gamma(\lambda_3+\lambda_4,\lambda_2-\lambda_3,
\lambda_2+\lambda_4+1,\lambda_1-\lambda_2).$
Similarly, the number of points in ${\bf N}^2$ that can be reached 
from $(\nu_2,\mu_2+1)$ inside $I_2\times I_5$
is given by
$\Gamma(\lambda_3+\lambda_4,\lambda_2-\lambda_3,
\lambda_2+\lambda_3+2,\lambda_1-\lambda_2).$
\end{proof}

\begin{corollary} \label{tres-dos}
Let 
$\mu=(\mu_1,\mu_2)$, $\nu=(\nu_1,\nu_2)$, and
$\lambda=(\lambda_1,\lambda_2)$  be partitions of $n$.
Assume that $\nu_2 \le \mu_2 \le \lambda_2$. 
Then
\[
\gamma^{\lambda}_{{\mu}{\nu}}=
(y-x)(y\ge x),
\]
where $x=\max \biggl(0, \biggl \lceil \frac {\mu_2+\nu_2+\lambda_2-n}{2}
 \biggl \rceil \biggl)$
and  $y=\biggl \lceil \frac {\mu_2+\nu_2-\lambda_2+1}{2} \biggl \rceil$.
\end{corollary}
\begin{proof}
Set $\lambda_3=\lambda_4=0$ in Theorem \ref{le:2rowst}. Then we
notice that the second possibility in the definition of $\Gamma$,
that is, when $c<y<c+d$, never occurs.
Note that $\nu_2+\mu_2-\lambda_2 \ge \nu_2+\mu_2-\lambda_1-1$ for all
partitions $\mu, \nu,$ and $\lambda$.
 Therefore, 
\[
\gamma^{\lambda}_{{\mu}{\nu}}=
 \sigma_{\lambda_2+1,1}(\nu_2+\mu_2-\lambda_2)
-\sigma_{\lambda_2+1,1}(\nu_2+\mu_2-\lambda_1-1).
\]

Suppose that $\nu_2+\mu_2-\lambda_2<0.$ From the definition of 
$\sigma_{\lambda_2+1,1}$ we obtain that $\gamma^{\lambda}_{{\mu}{\nu}}=0$.
Therefore, in order to have $\gamma^{\lambda}_{{\mu}{\nu}}$ not equal to 
zero, we should assume that $\nu_2+\mu_2-\lambda_2\ge 0.$

If $0\le \nu_2+\mu_2-\lambda_2 < \lambda_2+1$, then 
\[
\sigma_{\lambda_2+1,1}(\nu_2+\mu_2-\lambda_2)=\biggl \lceil 
\frac{\nu_2+\mu_2-\lambda_2+1}{2} \biggl \rceil
\]
Similarly, if  $0\le \nu_2+\mu_2-\lambda_2 < \lambda_2+1$, then
\[
\sigma_{\lambda_2+1,1}(\nu_2+\mu_2-\lambda_1-1)=\biggl \lceil 
\frac{\nu_2+\mu_2+\lambda_2-n}{2} \biggl \rceil
\]
It is easy to see that all other cases on the definition of $\sigma_{k,l}$ can not 
occur. Therefore, defining $x$ and $y$ as above, we obtain the desired result.
\end{proof}

\begin{Example}
 If  $\mu=\nu=\lambda=(l,l)$
or   $\mu=\nu=(2l,2l)$ and $\lambda=(3l,l)$,
then from the previous corollary, we obtain that
\begin{align*}
\gamma^{\lambda}_{{\mu}{\nu}}=\biggl\lceil \frac{l+1}{2} \biggl\rceil
                                  -\biggl\lceil \frac{l}{2} \biggl\rceil
                             =(l \text{  is even})
\end{align*}

Note that to apply Corollary \ref{tres-dos} to the second family of
shapes, we should first use the symmetries of the Kronecker product
to set $\nu=\lambda=(2l,2l)$ and $\mu=(3l,l)$.

\end{Example}

\begin{corollary} The Kronecker coefficients
$\gamma^{\lambda}_{{\mu}{\nu}}$, where $\mu$ and $\nu$ are two-row
partitions, are unbounded.
\end{corollary}
\begin{proof}
It is enough to construct an unbounded family of Kronecker
coefficients. Assume that $\mu=\nu=\lambda=(3l,l)$.
Then from the previous corollary we obtain that
\begin{align*}
\gamma^{\lambda}_{{\mu}{\nu}}=\biggl\lceil \frac{l+1}{2} \biggl\rceil
\end{align*}
\end{proof}

			\section{Sergeev's formula} \label{sergeev}

In this section we state Sergeev's formula for a Schur 
function of a difference of two alphabets. 
See \cite{Be-Ga}  or [7, section I.3] for
 proofs and  comments. 

\begin{definition}
Let $X_m = x_1+ \cdots+ x_m $ be a finite alphabet,
 and let \\$\delta_m=(m-1, m-2, \cdots, 1, 0)$. We
define 
$X_m^{\delta_m}$ by 
$X_m^{\delta_m}=x_1^{m-1}  \cdots x_{m-1}.$
\end{definition}

\begin{definition} 
Let $i(\alpha)$ denote the number of inversions of 
the permutation $\alpha$. We define the alternant 
to be
 \[
A_m^x P = \sum_{\alpha \in S_m}
(-1)^{i(\alpha)}
	P(x_{\alpha(1)},\cdots,x_{\alpha(m)}),
\]
 for any polynomial $P(x_1, \cdots, x_n)$. 
\end{definition}

\begin{definition}
Let  $\Delta$
be the operation of taking the 
 Vandermonde determinant of an alphabet, i.e.,
	\[
		\Delta(X_m)=\mathop{\rm det}(x_i^{m-j})_{i,j=1}^m.
	\]
\end{definition}

\begin{theorem}[Sergeev's Formula]
 Let $X_m = x_1+ \cdots + x_m , \; and \;\, Y_n=  y_1+
\cdots+ y_n $ be
	 two alphabets. Then
  		\[
			s_{\lambda}[X_m - Y_n] =
			\frac {1}
			{\Delta(X_m)\Delta(Y_n) }
			A_m^x A_n^yX_m^{\delta_m} Y_n^{\delta_n}
			\prod_{(i,j) \in \lambda} (x_i - y_j)
		\]
\end{theorem}

The notation $(i,j) \in \lambda $ means that the 
point $(i,j)$ belongs to the diagram of 
$\lambda$. We set $x_i=0$ for $i>m$ and $y_j=0$
 for $j>n$. 
	\medskip

We use Sergeev's formula as a tool for making some 
calculations we  need for the next two sections. 

\medskip

\begin{enumerate}

\item{}	Let ${\mu}=(1^{e_1} m_2)$ be a hook. (We are assuming that
$e_1 \ge 1$ and $ m_2 \ge 2$.) 
Let $X^1=\{x_1\}$ and	$X^2=\{x_2\}$.  
			\begin{equation}
   s_{\mu}[x_1 - x_2] = 
				(-1)^{e_1} x_1^{m_2-1} x_2^{e_1} (x_1 - x_2).
				\label{eq:hook}
			\end{equation}

\item{}
 Let ${\nu}=(\nu_1,\nu_2)$ be a two-row partition. 
Let $Y=\{y_1,y_2\}$.  Then
			\begin{equation}
				s_{\nu}[y_1+y_2] = 
    \frac{(y_1y_2)^{\nu_2}(y_1^{\nu_1-\nu_2+1}
	        -y_2^{\nu_1-\nu_2+1})}
    {y_1-y_2}.
				\label{eq:2row}
			\end{equation}

\medskip

\item{}	We say that a partition $\lambda$ is a double 
hook if $(2,2)\in\lambda$ and
it has the form ${\lambda}=(1^{d_1}2^{d_2}n_3 \,
n_4).$  In particular any two-row shape is considered
to be a double hook.

Let $\lambda$ be a double hook.
 Let $U=\{u_1,u_2\}$ and $V=\{v_1,v_2\}$. Then if $n_4 \ne 0$ then
$s_{\lambda}[u_1+u_2-v_1-v_2]$ equals
\begin{setlength}{\multlinegap}{40pt}
		\begin{multline}  
			\frac{(u_1-v_1)(u_2-v_1)(u_1-v_2)(u_2-v_2)}
  				{(u_1-u_2)(v_1-v_2)}(-1)^{d_1}
			(u_1u_2)^{n_3-2}(v_1v_2)^{d_2}\\
			\times
  		\left(u_2^{n_4-n_3+1}-u_1^{n_4-n_3+1}\right)
  		\left(v_2^{d_1+1}-v_1^{d_1+1}\right).
			\label{eq:double}
		\end{multline}
\end{setlength}
On the other hand, if $n_4=0$ then to
compute $s_{\lambda}[u_1+u_2-v_1-v_2]$
 we should write $\lambda$ as $(1^{d_1}2^{d_2-1}2\,n_3).$ That is, we set $d_1:=d_1$,
$d_2:=d_2-1$, $n_3:=2$, and $n_4=n_3$ in  (\ref{eq:double}). 

\medskip

\item{}	Let ${\lambda}$ be a hook shape, $\lambda=(1^{d_1}n_2). $
(We are assuming that $d_1 \ge 1$ and $ n_2 \ge 2$.) 
	Let $U=\{u_1,u_2\}$ and $V=\{v_1,v_2\}$.  Then 
$s_{\lambda}[u_1+u_2-v_1-v_2]$ equals
\begin{setlength}{\multlinegap}{40pt}
\begin{multline}
	(-1)^{d_1-1}\frac{1}{(u_1-u_2)}\frac{1}{(v_1-v_2)}
		\times\\
			\biggl\{ u_1v_1(u_1-v_1)(u_1-v_2)(u_2-v_1)
    u_1^{n_2-2}v_1^{d_1-1}
			 \\
  \;\;\;\;\;\;\;\;\;\;\;\;\;\;-u_1v_2(u_1-v_2)(u_1-v_1)(u_2-v_2)
   u_1^{n_2-2}v_2^{d_1-1}\\
  \;\;\;\;\;\;\;\;\;\;\;\;\;\;-u_2v_1(u_2-v_1)(u_2-v_2)
				(u_1-v_1)u_2^{n_2-2}v_1^{d_1-1}
			\\
	+ u_2v_2(u_2-v_2)(u_2-v_1)(u_1-v_2)
			u_2^{n_2-2}v_2^{d_1-1}.\biggr\}
			\label{eq:s-hook}
	 	\end{multline}
\end{setlength}
\end{enumerate}

	\section{The case of two hook shapes} \label{2hooks}

In this section we derive an explicit formula for the
 Kronecker coefficients $\gamma^{\lambda}_{{\mu}{\nu}}$ 
in the case in which
 ${\mu}=(1^e u)$, and ${\nu}=(1^f v)$
 are both hook shapes.  
Given a partition $\lambda$  the
 Kronecker coefficient $\gamma^{\lambda}_{{\mu}{\nu}}$
tells us whether  point $(u,v)$ belongs to some regions in ${\bf N}^2$
determined by ${\mu}$, ${\nu}$ and $\lambda$.

Recall that we denote the characteristic function 
by enclosing
a proposition $P$ with brackets, $(P)$.

\begin{lemma}\label{le:4}
Let $(u,v) \in {\bf N}^2$ and let $R$ be the rectangle with vertices
$(a,b)$, $(b,a)$, $(c,d)$, and $(d,c)$, with $a \ge b$, $c \ge d$,
$c \ge a$ and $d \ge b$. (Sometimes, when $c=d=e$, 
we denote this rectangle as
$(a,b;e)$.)

Then
 $(u,v) \in R $ if and only if
$|v-u|\le a-b $ and $ a+b \le u+v \le c+d$
\end{lemma}

\begin{proof}
Let $L_1$ be the line of slope $1$ passing through $(u,v)$, and let
$L_2$ be the line of slope $-1$ passing through $(u,v)$. Then we 
have that
\begin{align*}
&L_1 : y=x+v-u\\
&L_2 : y=-x+u+v
\end{align*}
The point $(u,v)$ is in $R$ if and only if
 $L_1$ is between the lines of slope $1$ passing throught
$(a,b)$ and $(b,a)$. That is, $a-b \le v-u \le b-a$ 
and  $L_2$ is between the lines of slope $-1$ passing throught
$(a,b)$ and $(c,d)$. That is, $a+b \le u+v \le c+d$
\end{proof}

\begin{theorem} \label{th:2h}
Let $\lambda$, $\mu$ and $\nu$ be partitions of $n$, where 
$\mu=(1^e u)$ and $\nu=(1^fv)$ are hook shapes.
 Then the Kronecker coefficients
$\gamma^{\lambda}_{{\mu}{\nu}}$  are given by the following:

\begin{enumerate}

\item If $\lambda$ is a one-row shape, then 
$\gamma^{\lambda}_{{\mu}{\nu}}=\delta_{\mu,\nu}.$

\item {} If $\lambda$ is not contained in a double hook
 shape, then 
$\gamma^{\lambda}_{{\mu}{\nu}}=0.$

\item {} Let $\lambda=(1^{d_1} 2^{d_2} n_3 n_4)$ be a double
hook. Let $x=2d_2+d_1$. Then
\begin{setlength}{\multlinegap}{50pt}
\begin{multline*}
\gamma^{\lambda}_{{\mu}{\nu}}=
(n_3 -1 \le \frac{e+f-x}{2} \le n_4) (|f-e| \le d_1)\\
+(n_3 \le \frac{e+f-x+1}{2} \le n_4) (|f-e| \le d_1+1).	
\end{multline*}
\end{setlength}
Note that if $n_4=0$, then we shall rewrite $\lambda=(1^{d_1}2^{d_2-1}2\,n_3)$
before using the previous formula.

\item {} Let $\lambda=(1^d w)$ be a hook shape.
Suppose that $e \le u$, $f \le v$, and $d \le w$. Then
\[
\gamma^{\lambda}_{{\mu}{\nu}}= (e \le d+f)(d \le e+f)( f \le e+d).
\]

\end{enumerate}
\end{theorem}

\begin{proof} 
Set $X=\{1,x\}$ and $Y=\{1,y\}$
 in  the comultiplication expansion (\ref{eq:main}) to
 obtain
\begin{equation}
  s_{\lambda}[(1-x)(1-y)] = 
		\sum_{\mu,\nu} \gamma^{\lambda}_{{\mu}{\nu}}
		s_{\mu}[1-x]s_{\nu}[1-y],
		\label{eq:main2hook}
	\end{equation}

We use equation (\ref{eq:hook}) to replace 
$s_{\mu}$ and $s_{\nu}$ 
in the right hand side of (\ref{eq:main2hook}). Then we
divide the resulting equation by $(1-x)(1-y)$ to get
	 \[
		\frac{s_{\lambda}[1-y-x+xy]}{(1-x)(1-y)} = 
		\sum_{\mu,\nu} \gamma^{\lambda}_{{\mu}{\nu}}
  (-x)^{e}(-y)^{f}. 
  \]
Therefore, 
\[
	\gamma^{\lambda}_{{\mu}{\nu}}=
\big[(-x)^{e}(-y)^{f}\big]\frac{s_{\lambda}[1-y-x+xy]}{(1-x)(1-y)},
\]
 when $\mu$ and $\nu$ are
 hook shapes.

{\bf Case 1.} If $\lambda$ is not contained in any double
hook, then the point $(3,3)$ is in $\lambda$, and  by
 Sergeev's formula, $s_{\lambda}[1-y-x+xy]$ equals
zero.

\medskip
{\bf Case 2.} Let $\lambda=(1^{d_1}2^{d_2}n_3\,n_4)$ be a
double hook.  
Set $u_1=1$, $u_2=xy$, $v_1=x$, and $v_2=y$ in 
(\ref{eq:double}). Then we divide by $(1-x)(1-y)$ on both
sides of the resulting equation to obtain
	\begin{multline}
		\frac{s_{\lambda}[1-y-x+xy]}{(1-x)(1-y)}
		=(-1)^{d_1} 
		(xy)^{n_3+d_2-1}\\
		\times(1-x)(1-y)
  \left( \frac{1-(xy)^{n_4-n_3+1}}{1-xy} \right)
		\left( \frac{x^{d_1+1}-y^{d_1+1}}{x-y} \right).
	\label{eq:case 2.1}
	\end{multline}
Note: If $n_4=0$ then we should write $\lambda=(1^{d_1}2^{d_2-1}2\,n_4)$
in order to use (\ref{eq:double}).

Let $p$ be a point in ${\bf N}^2.$ We say  that $(i,j)$ 
can be reached from $p$, written
$p \leadsto (i,j)$, if $(i,j)$ can be reached from $p$ by moving any
number of steps south-west or north-west.
We have defined a weight function by
\[
\omega_p(i,j)=
\begin{cases}
x^iy^j, &\text{ if } p \leadsto (i,j);\\
0, &\text{otherwise.}
\end{cases}
\]
Let 
$
\omega_p(T)=\sum_{(i,j) \in T}{\omega_p(i,j)}
$
be the generating function of a region $T$ in ${\bf N}^2$.
Let $R$ be the rectangle with vertices $(0,d_1)$, $(d_1,0)$,
$(d_1+n_4-n_3,n_4-n_3)$ and $(n_4-n_3,d_1+n_4-n_3)$.
Then 
\begin{align*}
\omega_{(d_1+n_4-n_3,n_4-n_3)}(R)&=
\left(
\frac{1-(xy)^{n_4-n_3+1}}{1-xy} 
\right)
\left( 
\frac{ x^{d_1+1}-y^{d_1+1}}{x-y}
\right)\\
&=\sum_{k=0}^{n_4-n_3}\sum_{i+j=d_1}{(xy)^k x^i y^j}. 
\end{align*}
 See Table 5.

\bigskip

\centerline{
\begin{tabular}{| c| c |c |c |c |c |c |c| c| } \hline
  &  &  &  & $1$&  &  &  &     \\  \hline
  &  &  & $1$&  & $1$&  &  &     \\   \hline
  &  & $1$&  & $1$&  & $1$&  &     \\   \hline
  & $1$&  & $1$&  & $1$&  & $1$&     \\  \hline
 $1$&  & $1$&  & $1$&  & $1$&  &$1$     \\ \hline
  & $1$&  & $1$&  & $1$&  & $1$&    \\ \hline
  &  & $1$&  & $1$&  & $1$&  &  \\ \hline
  &  &  & $1$&  & $1$&  & &     \\ \hline
  &  &  &  & $1$&  & &  &     \\ \hline
\end{tabular} }

\medskip
\centerline{ Table 5.  }
{\centerline {$d_1=4$ and $n_3-n_4=4$} }

\bigskip

Recall that we are using matrix coordinates, and that the
upper-left corner has coordinates $(0,0)$. The coordinates of the
four vertices of $R$ in Table $5$ are $(0,4)$, $(4,0)$, 
$(8,4)$, and $(4,8)$.

We interpret the right-hand side of
(\ref{eq:case 2.1})  as the sum of four different generating functions.
To be more precise, the right-hand side of (\ref{eq:case 2.1})
can be written as 
$
\sum_{i=1}^4 \omega_{p_i}(r_i)
$
where
$p_1=(n_4+d_2-1,n_4+d_2+d_1-1)$ and 
$R_1=\{n_3+d_2+d_1-1, n_3+d_2-1;n_4-n_3\}$,
$p_2=(n_4+d_2,n_34+d_2+d_1-1)$ and 
$R_2=\{n_3+d_2+d_1, n_3+d_2-1;n_4-n_3\}$,
$p_3=(n_4+d_2-1,n_4+d_2+d_1)$ and 
$R_3=\{n_3+d_2+d_1-1, n_3+d_2;n_4-n_3\}$,
and 
$p_4=(n_4+d_2+d_1,n_4+d_2+d_1)$ and 
$R_4=\{n_3+d_2+d_1, n_3+d_2+d_1;n_4-n_3\}$.

\smallskip

We observe that $R_1 \cup R_2$ (and $R_3 \cup R_4$ ) are rectangles in
${\bf N}^2$. Moreover,
\begin{equation}
\gamma_{\mu \nu}^{\lambda}=((e,f) \in R_1 \cup R_2)
+((e,f) \in R_3 \cup R_4).
\label{rect}
\end{equation}
The vertices of rectangle $R_1 \cup R_4$ are given (using the notation
of  
\ref{le:4}) by
\begin{align*}
a&=n_3+d_2+d_1-1
&b&=n_3+d_2-1\\
c&=n_4+d_2+d_1
&d&=n_4+d_2
\end{align*}
Similarly, the vertices of rectangle $R_2\cup R_3$ are given  by
\begin{align*}
&a=n_3+d_2+d_1
&b&=n_3+d_2-1\\
&c=n_4+d_2+d_1
&d&=n_4+d_2-1
\end{align*}

Applying Lemma \ref{le:4} to (\ref{rect}) we obtain 
\begin{setlength}{\multlinegap}{50pt}
\begin{multline*}
\gamma^{\lambda}_{{\mu}{\nu}}=
(n_3 -1 \le \frac{e+f-x}{2} \le n_4) (|f-e| \le d_1)\\
+(n_3 \le \frac{e+f-x+1}{2} \le n_4) (|f-e| \le d_1+1).	
\end{multline*}
\end{setlength}

{\bf  Case 3. $\lambda$ is a hook.}
Suppose that $\lambda$ is a hook, $\lambda=(1^{d} w)$.
 Set $u_1=1$, $u_2=xy$, $v_1=x$, and $v_2=y$ 
in (\ref{eq:s-hook}). Then we divide by $(1-x)(1-y)$ on
both sides of the resulting equation to obtain 
	\begin{multline}
 		\frac{s_{\lambda}[1-y-x+xy]}{(1-x)(1-y)}=
	 	(-1)^{d} \left( \frac{x^{d+1}-y^{d+1}}{x-y} \right) 
			\left( \frac{1-(xy)^{w}}{1-xy} \right)
			\\
 		+(-1)^{d-1} 
		xy \left( \frac{x^{d}-y^{d}}{x-y} \right)
 		\left( \frac{1-(xy)^{w-1}}{1-xy} \right).
	\label{eq:case2.2}
	\end{multline}

We want to interpret this equation as a generating function for a region $T$
using the weight $\omega$. We proceed as
follows:

Let $R_1$ be the rectangle with vertices $(d,0),$ $ (0,d),$ $(d+w-1,w-1)$, and
$(w-1,d+w-1)$. 
Then 
\begin{equation}
\omega_{(w-1,d+w-1)}(R_1)=
\left(
\frac{1-(xy)^{w}}{1-xy} 
\right)
\left( 
\frac{ x^{d+1}-y^{d+1}}{x-y}
\right)=\sum_{k=0}^{w-1}\sum_{i+j=d}{(xy)^kx^iy^j}.
\label{eq:serexp1} 
\end{equation}
(See Table 5.) Similarly, let $R_2$ be the rectangle with vertices 
$(d,1),$ $ (1,d),$ $(d+w-2,w-1)$, and
$(w-1,d+w-2)$. Then 
\begin{equation}
\omega_{(w-1,d+w-2)}(R_2)=xy
\left(
\frac{1-(xy)^{w-1}}{1-xy} 
\right)
\left( 
\frac{ x^{d}-y^{d}}{x-y}
\right)=xy\sum_{k=0}^{w-2}\sum_{i+j=d-1}{(xy)^kx^iy^j}.
\label{eq:serexp2} 
\end{equation}

Observe that the points in ${\bf N}^2$ that can be reached from $(0,d)$ in $R_1$
and the points in ${\bf N}^2$ that can be reached from $(1,d)$ in  $R_2$ are
disjoint. Moreover, they completely fill the rectangle $R_1 \cup R_2$.
See Table 6.

\bigskip

\centerline{
\begin{tabular}{ |c| c| c|c |c |c| c| c| c| c| }
\hline
  &  &  &  & $\phantom{+}1$& &  &  &  &    \\ \hline
  &  &  & $\phantom{+}1$&$-1$& $\phantom{+}1$
& &  & &    \\ \hline
  &  &
$\phantom{+}1$&$-1$&$\phantom{+}1$&$-1$&$\phantom{+}1$
& & &   \\ \hline
   &
$\phantom{+}1$&$-1$&$\phantom{+}1$&$-1$&$\phantom{+}1$
&$-1$&$\phantom{+}1$& &   \\ \hline
$\phantom{+}1$&$-1$&$\phantom{+}1$&$-1$&$\phantom{+}1$
&$-1$&$\phantom{+}1$&$-1$&$\phantom{+}1$&  \\ \hline
&$\phantom{+}1$&$-1$&$\phantom{+}1$&$-1$&
$\phantom{+}1$&$-1$&$\phantom{+}1$&$-1$&$\phantom{+}1$ 
\\ \hline
  & &$\phantom{+}1$&$-1$&$\phantom{+}1$&$-1$&
$\phantom{+}1$&$-1$&$\phantom{+}1$&   \\ \hline
  &  &
&$\phantom{+}1$&$-1$&$\phantom{+}1$&$-1$
&$\phantom{+}1$& 
&   \\ \hline
  &  &  &  &$\phantom{+}1$&$-1$&$\phantom{+}1$&  & 
&    \\  \hline
  &  &  &  &  & $\phantom{+}1$&  &  &  &    \\ \hline
\end{tabular}}

\bigskip
\centerline{Table 6}
\centerline{$d=4, w=6.$}
\bigskip

Note that $R_2$ is contained in $R_1$. We obtain that
\[
\omega_{(w-1,d+w-1)}(R_1)+\omega_{(w-1,d+w-2)}(R_2)=|(e,f) \in R_1|
\]
We use apply Lemma \ref{le:4} to the previous equation to obtain:
\[
(|e-f| \le d)(d \le e+f \le d+2w-2).
\]
But, by hypothesis, $e \le u$, $f \le v$, and $d \le w$. Therefore, 
 this system is equivalent to 
$
(d \le e+f)(f \le e+d) (e \le d+f)
$, as desired.
\end{proof}

\begin{corollary} \label{Co:hh2}
Let $\lambda$, $\mu$, and $\nu$ be partitions of $n$, where 
$\mu=(1^eu)$ and $\nu=(1^fv)$ are hook shapes and 
$\lambda=(\lambda_1,\lambda_2)$ is a two-row shape. Then
the Kronecker coefficients $\gamma^{\lambda}_{{\mu}{\nu}}$
are given by
\[
\gamma^{\lambda}_{{\mu}{\nu}}=(\lambda_2-1 \le e \le \lambda_1)(e=f)
+(\lambda_2 \le \frac{e+f+1}{2}  \le \lambda_1)(|e-f|\le 1).
\]
\end{corollary} 
\begin{proof}
In Theorem \ref{th:2h},  set $d_1=d_2=0$, $n_3=\lambda_2$ and $n_4=\lambda_1.$ 
\end{proof}
\begin{corollary}  Let $\lambda$, $\mu$ and $\nu$ be
partitions of
 $n$, where $\mu$ and $\nu$ are hook shapes. Then the
Kronecker coefficients are bounded. Moreover, the only
possible values for the Kronecker coefficients are $0$,
 $1$ or $2$. 
\end{corollary}

	\section{The case of a hook shape and a two-row 
					shape} ~\label{hook2row}

In this section we derive an explicit formula for
 the Kronecker coefficients in the case $\mu=(1^{e_1}m_2)$ is 
a hook and $\nu=(\nu_1,\nu_2)$ is a  two-row shape. Given a partition 
$\lambda$, the Kronecker coefficients 
$\gamma^{\lambda}_{{\mu}{\nu}}$ tell us whether the point 
$(e_1,\nu_2)$ belongs to some regions
in ${\bf N}^2$ determined by $\mu$, $\nu$ and $\lambda$.

Using the symmetry properties of the Kronecker product, we may
 assume that  if $\lambda=(1^{d_1}2^{d_2}n_3n_4)$ then 
 $n_4-n_3 \leq d_1$. (If $n_4=0$ then we should rewrite $\lambda$ as 
$(1^{d_1}2^{d_2-1}2\,n_3)$. Moreover, our hypothesis 
becomes $n_3-2 \leq d_1$.) 

Recall that we denote the value of the characteristic
function at proposition $P$ by $(P)$.

\begin{theorem}
 Let $\lambda$, $\mu$ and $\nu$ be partitions of $n$,
where
$\mu=(1^{e_1}m_2)$ is a hook and
$\nu=(\nu_1,\nu_2)$ is a two-row shape.  Then the Kronecker
coefficients
$\gamma^{\lambda}_{{\mu}{\nu}}$ are given by the 
following:

\begin{enumerate} 

\item If $\lambda$ is a one-row shape, then 
$\gamma^{\lambda}_{{\mu}{\nu}}=\delta_{\mu,\nu}.$

\item {} If $\lambda$ is not contained in any double hook,
 then $\gamma^{\lambda}_{{\mu}{\nu}}=0.$

\item {}  Suppose $\lambda=(1^{d_1}2^{d_2}n_3n_4)$ is
 a double hook.  Assume that $n_4-n_3 \leq d_1$.\\ 
(If $n_4=0$, then we should write
$\lambda=(1^{d_1}2^{d_2-1}2\,n_3)$.)  Then
\begin{setlength}{\multlinegap}{20pt}
\begin{multline}
\gamma_{\mu \nu}^{\lambda}=
(n_3 \leq \nu_2 -d_2 -1\leq n_4)(d_1+2d_2 < e_1 <d_1+2d_2+3)\\
+(n_3 \leq \nu_2-d_2 \leq n_4)(d_1+2d_2 \le e_1 \le d_1+2d_2+3)\\
+ (n_3 \leq \nu_2 -d_2 +1\leq n_4)(d_1+2d_2 < e_1 <d_1+2d_2+3)\\
-(n_3+d_2+d_1=\nu_2)(d_1+2d_2+1 \le e_1 \le d_1+2d_2+2).
\notag
\end{multline}
\end{setlength}

\item If $\lambda$ is a hook, see Corollary \ref{Co:hh2}.
\end{enumerate}
\end{theorem}
\begin{proof}  
 Set $X=1+x$ and $Y=1+y$ in the comultiplication
expansion (\ref{eq:main}) to obtain
\begin{equation}
s_{\lambda}[(1-x)(1+y)] = 
		\sum_{\mu,\nu} \gamma^{\lambda}_{{\mu}{\nu}}
		s_{\mu}[1-x]s_{\nu}[1+y].
		\label{eq:mainhookrow}
	\end{equation}

Use (\ref{eq:hook}) and
(\ref{eq:2row}) to replace $s_{\mu}$ and
$s_{\nu}$ in  the right-hand side of 
(\ref{eq:mainhookrow}), and divide by
 $(1-x)$  to obtain
	\begin{equation}
		\frac{s_{\lambda}[(1-x)(1+y)]}{1-x}
		= \sum_{\substack{\mu=(1^{e_1}m_2)\\ 
			\nu=(\nu_1, \nu_2)}}
			\gamma^{\lambda}_{{\mu}{\nu}}
			(-x)^{e_1} y^{\nu_2}
			\biggl(
			\frac{1-y^{\nu_1-\nu_2+1}}{1-y}
			\biggl). 
	\label{eq:sub2hooks}
	\end{equation}
If $\lambda$ is not contained in any double
hook, then the point $(3,3)$ is in $\lambda$, and  by
 Sergeev's formula, $s_{\lambda}[(1-x)(1+y)]$ 
equals zero.

Since we already computed the Kronecker coefficients when $\lambda$
is contained in a hook, we can assume 
for the rest of this proof that $\lambda$ is
a double hook. Let $\lambda=(1^{d_1}2^{d_2}n_3\,n_4)$.
(Note: If $n_4=0$ then we should write $\lambda=(1^{d_1}2^{d_2-1}2\,n_4).$)

 Set $u_1=1$, $u_2=y$, $v_1=x$,
 and $v_2=xy$ in (\ref{eq:double}), and multiply by
 $\frac{1-y}{1-x}$ on both sides of the resulting equation.
		\begin{setlength}{\multlinegap}{20pt}
		\begin{multline}
			\sum_{\substack{\mu=(1^{e_1}m_2)\\ 
			\nu=(\nu_1, \nu_2)}}
			\gamma^{\lambda}_{{\mu}{\nu}}
			(-x)^{e_1} y^{\nu_2}
			\big(
			1-y^{\nu_1-\nu_2+1}
			\big)
			=
			(y-x)(1-xy)(1-x)\\
			\times
			(-x)^{d_1+2d_2}y^{n_3+d_2-1}
			\left( \frac{(1-y^{n_4-n_3+1})(1-y^{d_1+1})}
			{1-y} \right).
  \label{eq:gen2}
		\end{multline}
		\end{setlength}
We have that $(y-x)(1-xy)(1-x)=y-x(1+y+y^2)+x^2(1+y+y^2)-x^3y.$ Therefore, 
looking at the coefficient of $x$ on both sides of the equation, we
see that $\gamma^{\lambda}_{{\mu}{\nu}}$ is zero if $e_1$
is different from $d_1+2d_2, d_1+2d_2+1, d_1+2d_2+2,$ or $d_1+2d_2+3.$

Let $e_1=d_1+2d_2$ or $e_1=d_1+2d_2+3$. Since $\nu_2 \le n/2$, we have that
		\begin{align*}
			\gamma_{\mu \nu}^{\lambda} 
			&=[y^{\nu_2}]
			\sum_{\substack{\mu=(1^{e_1}m_2)\\
			\nu=(\nu_1, \nu_2)}}
			\gamma^{\lambda}_{{\mu}{\nu}}
			y^{\nu_2}\\
			&=[y^{\nu_2}]
			\sum_{\substack{\mu=(1^{e_1}m_2)\\ 
			\nu=(\nu_1, \nu_2)}}
			\gamma^{\lambda}_{{\mu}{\nu}}
			y^{\nu_2}
			(1-y^{\nu_1-\nu_2+1}) 
			& & \text{($\nu_1+1 > n/2$)} \\
			&=
			[y^{\nu_2}]
			y^{n_3+d_2}(1-y^{d_1+1})
			\left( \frac{1-y^{n_4-n_3+1}}
			{1-y} \right)
                        & & \text{(Eq. \ref{eq:gen2})}\\
			&=[y^{\nu_2}]
			y^{n_3+d_2}(1-y^{d_1+1})
			\sum_{k=0}^{n_4-n_3}{y^k}
			\\
			&=[y^{\nu_2}]
			y^{n_3+d_2}
			\sum_{k=0}^{n_4-n_3}{y^k}.
			& & \text{($n_3+d_2+d_1 \ge n/2$)} 			
		\end{align*}
We have obtained that for $e_1=d_1+2d_2$ or $e_1=d_1+2d_2+3$
 \[
\gamma_{\mu \nu}^{\lambda}=
(n_3 \leq \nu_2-d_2 \leq n_4).
\]
Let $e_1=d_1+2d_2+1$ or $e_1=d_1+2d_2+2$. Since
$\nu_2 \le \lfloor \frac{n}{2} \rfloor$ we have that
{\allowdisplaybreaks
		\begin{align*}
			\gamma_{\mu \nu}^{\lambda}&=
			[y^{\nu_2}]
			\sum_{\substack{\mu=(1^{e_1}m_2)\\ 
			\nu=(\nu_1, \nu_2)}}
			\gamma^{\lambda}_{{\mu}{\nu}}
			y^{\nu_2}\\
			&=[y^{\nu_2}]
			\sum_{\substack{\mu=(1^{e_1}m_2)\\ 
			\nu=(\nu_1, \nu_2)}}
			\gamma^{\lambda}_{{\mu}{\nu}}(1-y^{\nu_1-\nu_2+1})
                        \\
			&=
			[y^{\nu_2}]
			y^{n_3+d_2-1}(1+y+y^2)(1-y^{d_1+1})
			\left( \frac{1-y^{n_4-n_3+1}}
			{1-y} \right)
                        \\
			&= \bigg(
			[y^{\nu_2}]
			y^{n_3+d_2-1}(1+y+y^2)
			\left( \frac{1-y^{n_4-n_3+1}}
			{1-y} \right) \biggl)
			-(n_3+d_2+d_1=\nu_2)\\
			&=\bigg( [y^{\nu_2}]
			y^{n_3+d_2-1}(1+y+y^2)
			\sum_{k=0}^{n_4-n_3}{y^k}
			\bigg) -(n_3+d_2+d_1=\nu_2)			
		\end{align*}
}
We have obtained that for $e_1=d_1+2d_2+1$ or $e_1=d_1+2d_2+2$
 \begin{multline}
\gamma_{\mu \nu}^{\lambda}=
(n_3 \leq \nu_2 -d_2 -1\leq n_4)
+(n_3 \leq \nu_2-d_2 \leq n_4)\\
+ (n_3 \leq \nu_2 -d_2 +1\leq n_4)-(n_3+d_2+d_1=\nu_2).
\end{multline}
\end{proof}

\begin{corollary}
The Kronecker cofficients, $\gamma^{\lambda}_{\mu \nu}$,
where $\mu$ is a hook and $\nu$ is a two-row shape are
always $0$, $1$, $2$
or $3$.
\end{corollary}

\section{Final comments}

The inner product of symmetric functions was discovered by J.
H. Redfield \cite{Red} in 1927, together with the scalar
product of symmetric functions. He called them cup and cap
products, respectively. D.E. Littlewood \cite{Li1, Li2} reinvented the inner product in 1956.

More recently, I.M. Gessel \cite{Ge} and A. Lascoux 
\cite{La} obtained combinatorial interpretations for  the
Kronecker coefficients in some restricted cases; Lascoux in
the case where $\mu$ and $\nu$ are hooks, and $\lambda$ a
straight tableaux, and Gessel in the case that $\mu$ and
$\nu$ are zigzag shapes and $\lambda$ is an arbitrary
 skew shape.
A. Lascoux interpreted the Kronecker coefficients, when
two of the shapes are hooks as counting clases of words
under some equivalence relation. We refer to \cite{La} or
\cite{Ga-Re} for a complete statement of his results. The
Corollary of Theorem 3 in this paper shows that each class
of words, under Lascoux's equivalence, contains either $0$,
$1$, or $2$ different representatives.
I. Gessel worked on a more general framework, contemplating
the occurrence of skew tableaux. It was shown in
\cite{Ga-Re}  that in the case where two of the partitions
are hook shapes, and the third one is an arbitrary straight
shape, his result is equivalent to Lascoux's.

 In \cite{Ga-Re}, A.M. Garsia and J.B. Remmel founded a
way to relate shuffles of permutations and Kronecker 
coefficients. From here they obtained a combinatorial
interpretation for the Kronecker coefficients when
$\lambda$ is a product of homogeneous  symmetric functions,
and $\mu$ and $\nu$ are arbitrary skew shapes. They also
showed how Gessel's and Lascoux's results are related.

J.B. Remmel \cite{Re1}, \cite{Re2}, and J.B. Remmel and T.
 Whitehead \cite{Re-Wh}, obtained formulas for  computing the
Kronecker coefficients in the same cases considered  in
this paper.  Their approach was mainly combinatorial:

First, they expanded the Kronecker product 
$s_{\mu}*s_{\nu}$ in terms of Schur functions using 
the Garsia-Remmel algorithm \cite{Ga-Re}. The problem of
computing the Kronecker coefficients was reduced to computing
signed sums of certain products of skew Schur functions.

Then they obtained a description of the coefficients 
that arise in the expansion of the resulting product of skew
Schur functions in terms of counting {\em 3-colored diagrams}
in \cite{Re1}, and
\cite{Re2} or {\em 4-colored diagrams} in \cite{Re-Wh}. At 
this point, they reduced the problem to computing a signed sum
of {\em colored diagram.}

Finally, they defined involutions on these signed sums
 to cancel negative terms, and obtained the desired formulas
by counting classes of restricted {\em colored diagram.} 

In general, it is not obvious how to go from the determination of the 
Kronecker coefficients $\gamma_{\mu,\nu}^{\lambda}$ when 
$\mu$ and $\nu$ are two-row shapes found in this paper,
and the one obtained by  J.B. Remmel and T. Whitehead \cite{Re-Wh}.
But, in some particular cases this is easy to see. For instance,
 when $\lambda$ is also a two-row shape, both formulas are exactly the same.


\begin{thebibliography}{99}

		\bibitem{Be-Ga} N. Bergeron and A.M. Garsia, ``Sergeev's		
		Formula and the Littlewood-Richardson Rule," {\em Linear and 	
		Multilinear Algebra} {\bf27}, 1990, pp. 79--100.

		\bibitem{Ga-Re} A.M. Garsia and J.B. Remmel, ``Shuffles 
		of permutations and Kronecker products," {\em Graphs Combin.} 
		1985, pp. 217--263.
   
  		\bibitem{Ge} I.M. Gessel, ``Multipartite P-partitions
		 and inner products of Schur functions," {\em Contemp. Math.}
		1984, pp. 289--302.

 		 \bibitem{La} A. Lascoux, ``Produit de Kronecker des 
		representations du group symmetrique,'' 
		{\em Lecture Notes in Mathematics } 1980, {\bf 795}, Springer Verlag
		pp. 319--329.

  		\bibitem{Li1} D.E. Littlewood, ``The Kronecker product 
		of symmetric group representations," {\em J. London Math. Soc.}
		 {\bf 31}, 1956, pp. 89--93.

  		\bibitem{Li2} D.E. Littlewood, ``Plethysm and inner 
		product of S-functions," {\em J. London Math. Soc.} {\bf 32},
		1957, pp. 18--22.

		\bibitem{Mac} I.G. Macdonald, {\em Symmetric Functions and 
		Hall Polynomials,} second edition, Oxford University Press, 1995.

  		\bibitem{Red} J.H. Redfield ``The theory of group reduced
		distribution," {\em Amer. J. Math.} {\bf 49}, 1927, pp. 433-455.

		\bibitem{Re1} J.B. Remmel, ``A formula for the Kronecker
		product of Schur functions of hook shapes," {\em J. Algebra} {\bf120}, 
		1989, pp. 100--118.

		\bibitem{Re2} J.B. Remmel,  ``Formulas for the expansion of 
		the Kronecker products $S_{(m,n)} \otimes S_{(1^{p-r},r)}$ 
		and $S_{(1^k 2^l)}\otimes S_{(1^{p-r},r)}$," {\em Discrete Math.} 
		{\bf99}, 1992, pp. 265--287.

		\bibitem{Re-Wh} J.B. Remmel and T. Whitehead,  ``On the 
		Kronecker product of Schur functions of two row shapes," 
		{\em Bull. Belg. Math. Soc. Simon Stevin} {\bf1}, 1994, pp. 649--683.

		\bibitem{Sag} B.E. Sagan, {\em The Symmetric Group}, Wadsworth 
		\& Brooks/Cole, Pacific Grove, California, 1991.
	\end{thebibliography}
\end{document}